\documentclass[11pt,leqno]{article}
\usepackage{amsthm,amsfonts,amssymb,amsmath,oldgerm}
\usepackage{epsfig}
\numberwithin{equation}{section}
\usepackage[thinlines]{easybmat}

\def\eps{\varepsilon }

\def\eps{\varepsilon}

\newcommand\br{\begin{remark}}
\newcommand\er{\end{remark}}
\newcommand\bp{\begin{pmatrix}}
\newcommand\ep{\end{pmatrix}}
\newcommand\be{\begin{equation}}
\newcommand\ee{\end{equation}}
\newcommand\ba{\begin{equation}\begin{aligned}}
\newcommand\ea{\end{aligned}\end{equation}}

\newcommand{\bap}{\begin{app}}
\newcommand{\eap}{\end{app}}
\newcommand{\begs}{\begin{exams}}
\newcommand{\eegs}{\end{exams}}
\newcommand{\beg}{\begin{example}}
\newcommand{\eeg}{\end{exaplem}}
\newcommand{\bpr}{\begin{proposition}}
\newcommand{\epr}{\end{proposition}}
\newcommand{\bt}{\begin{theorem}}
\newcommand{\et}{\end{theorem}}
\newcommand{\bc}{\begin{corollary}}
\newcommand{\ec}{\end{corollary}}
\newcommand{\bl}{\begin{lemma}}
\newcommand{\el}{\end{lemma}}
\newcommand{\bd}{\begin{definition}}
\newcommand{\ed}{\end{definition}}
\newcommand{\brs}{\begin{remarks}}
\newcommand{\ers}{\end{remarks}}

\newtheorem{theo}{Theorem}[section]

\newtheorem{exams}[theo]{Examples}

\numberwithin{equation}{section}

\newcommand{\D }{\mathcal{D}}
\newcommand{\U }{\mathcal{U}}
\newcommand{\A }{\mathcal{A}}

\newcommand{\CC}{{\mathbb C}}

\newtheorem{theorem}{Theorem}[section]
\newtheorem{proposition}[theorem]{Proposition}
\newtheorem{corollary}[theorem]{Corollary}
\newtheorem{lemma}[theorem]{Lemma}
\newtheorem{definition}[theorem]{Definition}

\newtheorem{example}[theorem]{Example}
\newtheorem{remark}[theorem]{Remark}

\newcommand\cA{{\cal  A}}
\newcommand\cB{{\cal  B}}
\newcommand\cD{{\cal  D}}
\newcommand\cU{{\cal  U}}
\newcommand\cH{{\cal  H}}

\newcommand\cK{{\cal  K}}

\newcommand\cF{{\cal  F}}

\newcommand\cM{{\mathcal M}}

\newcommand{\RM}{\mathbb{R}}
\newcommand{\ZM}{\mathbb{Z}}

\newcommand{\CM}{\mathbb{C}}

\newcommand{\NM}{\,\mathbb{N}}

\newcommand{\tr}{\,\mbox{\rm tr}}

\newcommand{\lb}{\label}

\newcommand{\bi}{\bibitem}

\newcommand{\beq}{\begin{equation}}
\newcommand{\eeq}{\end{equation}}

\title{
Convergence of Hill's method for nonselfadjoint operators
}

\author{\sc \small
Mathew A. Johnson\thanks{Indiana University, Bloomington, IN 47405;
matjohn@indiana.edu: Research of M.J. was partially supported by an NSF Postdoctoral Fellowship under NSF grant DMS-0902192.}
~~~~~
Kevin Zumbrun\thanks{Indiana University, Bloomington, IN 47405;
kzumbrun@indiana.edu: Research of K.Z. was partially supported
under NSF grants no. DMS-0300487 and DMS-0801745.
 }}
\begin{document}

\maketitle


\begin{center}
\end{center}


\begin{abstract}
By the introduction of a generalized Evans function defined
by an appropriate $2$-modified Fredholm determinant,
we give a simple proof of convergence in location and multiplicity
of Hill's method for numerical approximation of
spectra of periodic-coefficient ordinary differential operators.
Our results apply to operators of nondegenerate type, under the
condition that the principal coefficient
matrix be symmetric positive definite
(automatically satisfied in the scalar case).
Notably, this includes a large class of nonselfadjoint operators,
which were previously not treated.
The case of general coefficients depends on an interesting
operator-theoretic question regarding properties of Toeplitz matrices.
\end{abstract}

\section{Introduction}
The study of stability of spatially periodic traveling wave
solutions to various classes of partial differential equations
motivates the study of $L^2(\RM;\CC^n)$ (essential) spectra of
periodic-coefficient differential operators
\be\label{e:L1}
L=(\partial_x)^m a_m(x) + \dots + \partial_x a_1(x) + a_0(x)
\ee
on the line, where coefficients $a_j\in \CC^{n\times n}$ are periodic
with period $X$.
By Floquet theory,
it is equivalent to study the $L^2([0,X]; \CC^n)$ point spectra of
the family of
Bloch operators
$$
 L_\sigma=(\partial_x+i\sigma)^ma_m(x) + \dots
+(\partial_x+i\sigma)a_1(x) + a_0(x),
$$
where $X$ is the common period of the coefficients and $\sigma\in [0,2\pi)$ acts as a parameter.
Indeed, using this decomposition we have\footnote{Unless otherwise stated, throughout this paper all
functions are assumed to be complex valued and we adopt the notation $L^2(\RM)=L^2(\RM;\CM)$
and similarly for $L^2_{\rm per}([0,X])$.}
\[
{\rm spec}_{L^2(\RM)}\left(L\right)=\bigcup_{\sigma\in[0,2\pi)}{\rm spec}_{L^2_{\rm per}([0,X])}(L_\sigma);
 \]
see, for example, \cite{G} for more details.

Due to the mathematical difficulties
involved
in
analytically computing the $L^2(\RM)$ spectrum of such an,
in general, variable-coefficient and vector-valued, operator,
or, equivalently,
computing the periodic spectra of the full family of associated Bloch operators, 
the determination of spectrum of periodic-coefficient operators
is typically carried out numerically.
This may be accomplished in a number of ways: for example, shooting,
discretization, or various spectral and Galerkin methods.
See Appendix B, \cite{JZN}, for
further discussion.

A particularly natural
and direct
approach is Hill's method
\cite{DK},\footnote{
A convenient implementation may be found in the
numerical package SpectrUW \cite{CDKK}.}
%
a spectral Galerkin method carried out in a periodic Fourier basis, which
is exact in the constant-coefficient case.
In this method, to
approximate the spectra of $L_\sigma$ for a fixed $\sigma\in[0,2\pi)$,
one considers the eigenvalue problem
\begin{equation}\label{e:gspec}
L_\sigma v=\lambda v,
\end{equation}
by expressing the coefficients $a_j$ of $L_\sigma$ and the function $v$ as Fourier series in $L^2_{\rm per}([0,X])$,
as an infinite-dimensional matrix equation in $\ell^2$.
Truncating the Fourier modes to frequencies $|k|\leq J$ for each $J\in\NM$,
one then obtains a sequence of finite-dimensional
matrix eigenvalue problem whose eigenvalues approximate true eigenvalues of the operator $L_\sigma$ on $L^2_{\rm per}([0,X])$.
%
%
See Section \ref{s:hill} for further details.


This method is fast and easy to use, and in practice
appears to give excellent results
under quite general circumstances \cite{DK,BJNRZ1}.
However, up to now, an accompanying rigorous convergence theory has
been established only in certain commonly occurring
but restricted cases \cite{CuD}.
By convergence, we mean roughly that not
only is Hill's method accurate, meaning that the numerically computed eigenvalues are always close to the
actual eigenvalues of the associated Bloch-operator
(the ``no-spurious modes condition" of \cite{CuD}),
but also that the method is complete in the sense
that it faithfully produces all of $\sigma(L_\sigma)$ for a fixed $\sigma$: see \cite{CuD} for
a more precise discussion of convergence from this point of view.
Here, we make the simpler,
operational definition that on any bounded domain
$B=\{\lambda:\, |\lambda|\le R\}$
whose boundary contains no eigenvalue of $L_\sigma$,
the set of approximate eigenvalues lying in $B$ converges to the
set of exact eigenvalues of $L$ in both location and number;
see Cor.  \ref{Lconvthm}.\footnote{This
includes and slightly strengthens the definition of \cite{CuD}.}

Despite its obvious practical interest,
up to now the convergence of Hill's method
has been established to our knowledge
only for self-adjoint operators with principal coefficient $a_m=I$
\cite{CuD}.  In particular, though accuracy of Hill's method was
shown in \cite{CuD} under quite general assumptions,
completeness of the method in the non-selfadjoint
case, which arises naturally, for example, in the applications in \cite{BJNRZ1,BJNRZ2},
does not seem to have been fully addressed.

%

In this short paper,
we give a brief and simple proof of the convergence
of Hill's method applying 
to the general class of operators \eqref{e:L1}
such that $a_m$ is symmetric positive definite.
In the scalar case, this condition 
on the principal coefficient $a_m$ amounts to the mild requirement
that the operator be nondegenerate type.
In the system case, it is a
genuine restriction,
and it is an
interesting and apparently nontrivial
question,
related to certain properties of Toeplitz matrices,
to what extent the condition can be relaxed.
Notably, our
analysis
applies to the important case where the operator $L_\sigma$ is non-selfadjoint.

The main ingredient of our
our proof is the introduction of a generalized periodic Evans
function, of interest in its own right, consisting of a
$2$-modified Fredholm determinant $D_\sigma$ of an associated
Birman--Schwinger type operator, whose roots we show to agree in
location and multiplicity with the eigenvalues of $L_\sigma$.
For related analysis in the solitary wave case, see \cite{GLZ}.
Once these properties are established, the desired convergence
follows 
immediately by the observation that the corresponding 2-modified characteristic polynomial of the
$J^{\textrm{th}}$ Galerkin-truncation of $(L_\sigma-\lambda)v=0$ are
a subclass of the
approximants used to define the aforementioned $2$-modified Fredholm determinant
in the limit as $J\to \infty$, and furthermore that these approximates are a sequence of 
analytic functions converging locally uniformly to the generalized periodic Evans function.

A novel feature of the present analysis is that our argument yields
convergence of the spectrum in both location and multiplicity,
whereas the results of \cite{CuD} concerned only location.
On the other hand, there was established in \cite{CuD} a fast rate of
convergence to the smallest (in modulus)  eigenvalue in the self-adjoint
case, whereas our methods do not 
readily
appear to yield a rate.
A second novelty of our work is to make the connection to the Evans function,
putting the work in a broader context.


%

\section[Hilbert-Schmidt operators]{Hilbert--Schmidt operators and $2$-modified Fredholm determinants}
We begin
by recalling the basic properties
of $2$-modified Fredholm determinants, defined for
Hilbert--Schmidt perturbations of the identity;
see \cite{GGK1,GGK2}, \cite[Ch.\ XIII]{GGK3},
\cite[Sect.\ IV.2]{GK}, \cite{Si1}, \cite[Ch.\ 3]{Si2}
\cite[Sect.\ 2]{GLZ} for more details.

For a given Hilbert space $\cH$,\footnote{Throughout this paper, we will always assume that our Hilbert spaces are separable.}
 the Hilbert--Schmidt class
$\cB_2(\cH)$ is defined as the set of all bounded linear operators $A$ on $\cH$
for which the norm
\[
\|A\|_{\cB_2(\cH)} :=\sum_{j,k} |\langle Ae_j,e_k\rangle|^2=\tr_\cH ( A^* A)
\]
is finite, where $\{e_j\}$ is any orthonormal basis.
Evidently, $\|\cdot\|_{\cB_2(\cH)}$ is independent of the basis chosen.
Moreover, every operator in $\cB_2(\cH)$ is compact (Fredholm).

On a finite-dimensional
space $\cH$, we define the $2$-modified Fredholm determinant as
\ba \lb{2.34}
 {\det}_{2,\cH} (I_{\cH}-A):= {\det}_{\cH}((I_{\cH}-A)e^{A})
={\det}_{\cH}(I_{\cH}-A) \, e^{\tr_{\cH}(A)},
\ea
where $\det_{\cH}$
and $\tr_\cH$
denotes the usual determinant
and trace, respectively.
From this definition,
we have the useful estimates
	\be \lb{2.34b}
 |{\det}_{2,\cH}(I_{\cH}-A)| \leq e^{C\|A\|_{\cB_2(\cH)}^2}
\ee
and
\be\label{compare}
 |{\det}_{2,\cH}(I_{\cH}-A) - {\det}_{2,\cH}(I_{\cH}-B)| \leq \|A-B\|_{\cB_2(\cH)}
e^{C[\|A\|_{\cB_2(\cH)}+\|B\|_{\cB_2(\cH)}+1]^2},
\ee
where $C>0$ is a constant independent of the dimension of $\cH$.

To extend this notion of a determinant to an infinite dimensional Hilbert space $\cH$, we note that
for any $A\in\cB_2(\cH)$ the estimate \eqref{compare} allows us to
define the
%
$2$-modified Fredholm determinant unambiguously 
as the limit
\be\label{limdef}
 {\det}_{2,\cH}(I_{\cH}-A):= \lim_{J\to \infty} {\det}_{2,\cH_J}(I_{\cH_J}-A_J),
\ee
where $\cH_J$ is any increasing sequence of finite-dimensional subspaces
filling up $\cH$, and $A_J$ denotes the Galerkin approximation
$P_{\cH_J}A|_{\cH_J}$, where $P_J:\cH\to\cH_J$ is the orthogonal projection
onto $\cH_J$.
That is, thinking of the infinite-dimensional matrix representation
of $A$, the 2-modified Fredholm determinant is defined as the limit of such determinants
on finite, $J$-dimensional, minors as $J\to \infty $.

Alternatively, denoting the (countably many, since $A$ is Fredholm) eigenvalues
of $A$ as $\{\alpha_j\}_{j=1}^\infty$, and taking $\cH_J$ to be the
(total) eigenspace
associated with the eigenvalues $\{\alpha_j\}_{j=1}^J$ 
we find that
\be\label{productformula}
 {\det}_{2,\cH}(I_{\cH}-A)=\lim_{J\to\infty}\prod_{k=1}^J(1-\alpha_k)e^{\alpha_k},
\ee
which,
by
$\Pi_k(1-\alpha_k)e^{\alpha_k}\lesssim
\Pi_k (1+\alpha_k^2)
\sim e^{\sum_k \alpha_k^2}\le
e^{\|A\|_{\cB_2(\cH)}} $,
 is readily seen to converge for all $A\in\cB_2(\cH)$
by Weyl's inequality
$\sum |\alpha_j|^r\le \sum|s_j|^r$ for $r\ge 0$, where
$s_j$ denote the eigenvalues of $|A|:=(A^*A)^{1/2}$ \cite{Si1,W}.
%
This shows how the renormalization of the
standard determinant
$\det(I_{\cH}-A):=\Pi_j(1- \alpha_j)$
by factor $e^{\tr_\cH(A)}$ cancels the possibly divergent first-order terms
in $ \Pi_k (1-\alpha_k) \sim e^{\sum_k \alpha_k}$,
allowing the treatment of operators $A$
that are not in trace class $\cB_1:=\{A:\,
\||A|^{1/2}\|_{\cB_2(\cH)} <+\infty\}$.\footnote{
For $A\in \cB_1$, $ \tr_\cH(A)=\sum_j \alpha_j$ is absolutely convergent,
by Weyl's inequality with $r=1$, and so
the standard determinant
$\det_\cH (I_{\cH}-A)= \Pi_j (1-\alpha_j)$ converges.
For $A$ self-adjoint,
$\|A\|_{\cB_1}:=\||A|^{1/2}\|_{\cB_2(\cH)}=\sum_{j}|\alpha_j|$
and $\|A\|_{\cB_2(\cH)}=\sum_j|\alpha_j|^2$.
}

%

\bpr\label{property}
For $A\in \cB_2(\cH)$, the operator $(I_\cH-A)$ is invertible if and only if
${\det}_{2,\cH}(I_\cH-A)$ is non-zero.
\epr

\begin{proof}
By standard Fredholm theory, this is equivalent to the statement that
$0$ is an eigenvalue of $(I_\cH-A)$ if and only if $\det_{2,\cH}(I_\cH-A)= 0$.
Note that, since $A$ is Fredholm, it possesses a countable number of
isolated eigenvalues $\{\alpha_j\}$ of finite multiplicity, except
possibly at zero.
Choosing $J\in\mathbb{N}$ sufficiently large, then,
we may factor the product formula \eqref{productformula} as
\[
{\det}_{2,\cH}(I_\cH-K)=\left(\prod_{j=1}^J (1-\alpha_j)e^{\alpha_j}\right)
\left(\prod_{j=J+1}^\infty (1-\alpha_j)e^{\alpha_j}\right),
\]
where
\[
\prod_{j=J+1}^\infty(1-\alpha_j)e^{\alpha_j}\approx e^{\sum_{j=J+1}^\infty\alpha_j^2}\neq 0.
\]
It follows then that 
${\det}_{2,\cH}(I_\cH-A)$ vanishes if and only if $1-\alpha_j=0$
for some $1\le j\le J$, hence, since $J\in\mathbb{N}$ was arbitrary, if and only if $0$ is an eigenvalue of $(I_\cH-A)$.
\end{proof}

\section{Analysis of a simple case}\label{s:simple}

With the above preliminaries in
hand,
we now turn to our proof of convergence.
As a first step in this analysis, we
present a complete proof
in the
case of a second-order
operator with identity principal
part.  In later sections, we will then describe the extension of this proof to more
general cases, noting that most of the ideas can be found in this simpler context.

Consider a periodic-coefficient differential operator
\[
L_\sigma=(\partial_x+i\sigma)^2 + (\partial_x+i\sigma)a_1(x)
+ a_0(x)
\]
acting on vector-valued functions in
$L^2_{\rm per}([0,X])$,
$\sigma \in [0,2\pi)$ the Floquet parameter and
$a_j\in L^2([0,X])$ matrix-valued and periodic on $x\in [0,X]$.
We can rewrite this more generally as a family of operators in the simpler form
\be\label{e:L}
L_\sigma= \partial_x^2 + \partial_x A_1(\sigma,x) + A_0(\sigma,x),
\ee
where
$$
A_1= a_1 + 2i\sigma,\quad A_0= a_0-\sigma^2 + i\sigma a_1.
$$
In order to analyze the (necessarily discrete) spectrum of the operator $L_\sigma$, we
introduce a generalization of the periodic Evans function, a complex analytic function whose
roots coincide in location and multiplicity with the eigenvalues of $L_\sigma$ \cite{G},
expressed in terms of a 2-modified Fredholm determinant.  To
this end, notice that associated with the eigenvalue problem
\be\label{eig}
(L_\sigma-\lambda)U=0
\ee
is the equivalent problem
\be\label{equiv}
 (I+K(\sigma, \lambda))U=0,
\ee
where here $I$ is the identity operator on
$L^2_{\rm per}([0,X])$
and $K=K_1+K_0$, with
\[
K_1=\partial_x (\partial_x^2-1)^{-1} A_1, \quad
K_0= (\partial_x^2-1)^{-1} (A_0+1 -\lambda).
\]
In particular, notice that $\lambda$ is an eigenvalue of $L_\sigma$ if and only if $0$ is an eigenvalue of the operator $(I+K(\sigma,\lambda))$.
Before we can define the appropriate generalization of the Evans function, we need the following fundamental lemma.

\bl\label{l:hs}
For
$A_j\in L^2_{\rm per}([0,X])$,
the operator $K$ is Hilbert-Schmidt.
\el

\begin{proof}
Expressing $K_m$ in matrix form $\cK_m$ with respect to
the infinite-dimensional Fourier basis, we find that the corresponding matrix elements
can be expressed as
$$
[\cK_{1}]_{j,k}= \frac{ij}{1+j^2} \hat A_{1}(j-k),
$$
where $\hat A_1(m)$ denotes the $m^{th}$ Fourier coefficient of $A_1$, and $i:=\sqrt{-1}$.
Computing explicitly, we find by Parseval's Theorem that\footnote{Henceforth, Hilbert-Schmidt spaces $\cB_2$
will always be considered on the Hilbert space $L^2_{\rm per}([0,X])$.  That is, we adopt the notation
$\cB_2:=\cB_2(L^2_{\rm per}([0,X]))$.}
$$
\begin{aligned}
\|K_1\|_{\cB_2}&=\|\cK_1\|_{\cB_2}=
\sum_j
\frac{j^2}{(1+j^2)^2} \sum_k |\hat A_1(j-k)|^2\\
& =
\sum_j\frac{j^2}{(1+j^2)^2} \|A_1\|_{L^2_{\rm per}([0,X])}
<+\infty,
\end{aligned}
$$
hence $K_1$ is a Hilbert-Schmidt operator.  Similarly, we find that $K_0$ is Hilbert--Schmidt,
with norm
\[
\|K_0\|_{\cB_2}=\sum_{j}\frac{1}{(1+j^2)^2}\sum_k\left|\hat A_0(j-k)+(1-\lambda)\delta_j^k\right|^2,
\]
which implies that
$K=K_1+K_0\in \cB_2$
as claimed.
\end{proof}

\br
{
On the other hand, $K_1$ is not trace class if $\hat A_{1}(0):=\int_0^X A_1(x)dx\ne0$,
since then
$\sum_j |\cK_{1,jj}|= |\hat A_{1}(0)| \sum_j \frac{|j|}{1+ |j|^{2}} =+\infty$.
This illustrates the necessity of our extension of the usual notion of a determinant
to operators in $\cB_2$.
}
\er

\subsection{Generalized Periodic Evans Function}

By Lemma \ref{l:hs} in conjunction with Proposition \ref{property}, it follows that the zero eigenvalues
of $(I_{L^2_{\rm per}([0,X])}-K(\sigma,\lambda))$ can be identified through the use of a 2-modified Fredholm determinant.  This
leads us to the following definition.

\begin{definition}
For a fixed $\sigma\in[0,2\pi)$, we define the generalized periodic Evans function $D_\sigma:\CM\to\CM$ by
\be\label{evans2}
D_\sigma(\lambda):={\det}_{2,L^2_{\rm per}([0,X])}(I_{L^2_{\rm per}([0,X])}-K( \sigma, \lambda)).
\ee
\end{definition}

For ease of notation, throughout the rest of our analysis we will drop the dependence on the Hilbert space $L^2_{\rm per}([0,X])$
on the identity operator and all $2$-modified Fredholm determinants.  In particular, we will write
$D_\sigma(\lambda)=\det_2(I-K(\sigma,\lambda))$ for the above generalized Evans function.


\bt\label{evansthm}
For $A_j\in L^2_{\rm per}([0,X])$, the function 
$D_\sigma$ is complex-analytic in $\lambda$ and continuous in the parameter $\sigma$. 
Furthermore, the
roots of 
$D_\sigma$ for a fixed $\sigma\in[0,2\pi)$ correspond in location and multiplicity
with the eigenvalues of $L_\sigma$.
\et

\begin{proof}
Following the notation in Lemma \ref{l:hs}, for each $J\in\mathbb{N}$ we let $\cK_J:=([\cK]_{j,k})_{|j|,|k|\leq J}$
be the finite dimensional Galerkin matrix approximation
of the bi-infinite dimensional matrix representation of the operator $K$ defined above.
Clearly, then, for each fixed 
$J\in\mathbb{N}$ the finite-dimensional
approximation
$\Delta_J(\sigma,\lambda):=\det_{2}(I-\cK_J(\sigma,\lambda))$
is complex-analytic
in $\lambda$ and continuous in $\sigma\in[0,2\pi)$.  Furthermore, as in the proof of Lemma \ref{l:hs}
we have
\begin{equation}\label{e:kconverge}
\|\cK_{1,J}(\sigma,\lambda)-\cK_1(\sigma,\lambda)\|_{\cB_2}\leq\|A_1\|_{L^2([0,X])}\sum_{|j|\geq J+1}\frac{j^2}{(1+j^2)^2},
\end{equation}
where $\cK_{1,J}$ denotes the truncation of $\cK_1$,
and hence we find that $\cK_{1,J}\to\cK_1$ in
$\cB_2$
uniformly in both $\sigma$ and $\lambda$.
Similarly, we find that $\cK_{0,J}(\sigma,\lambda)\to\cK_0(\sigma,\lambda)$ in
$\cB_2$
uniformly in $\sigma$ and locally uniformly in $\lambda$,
and hence
the estimate \eqref{compare} implies\footnote{To use the estimate \eqref{compare} directly, one should
consider the operator $\cK_J$, which is technically defined on the finite-dimensional subspace $H_J$,
as being defined on the larger space $L^2_{\rm per}([0,X])$.  Throughout the remainder of our analysis
we will consider this extension without reserve.}
that $\Delta_J\to D_\sigma$
locally uniformly in $\lambda\in \CM$ and uniformly in $\sigma\in[0,2\pi)$.
It follows that the function $(\sigma,\lambda)\mapsto D_\sigma(\lambda)$
inherits the same regularity properties in $\lambda$ and $\sigma$ as the limiting
sequence $\Delta_J$, thus verifying the first claim of the Theorem.


Next, by equivalence of the problems \eqref{eig} and \eqref{equiv} together
with Proposition \ref{property}, we obtain immediately correspondence
in location of the roots of $D_\sigma$ and the eigenvalues of the operator $L_\sigma$.
To obtain agreement in multiplicity, consider an eigenvalue $\lambda_*$
of $L_\sigma$, with corresponding eigenspace $H_*$. Recalling that, by standard Fredholm theory, the eigenvalues
of $L_\sigma$ are countable, isolated, and have finite-multiplicity\footnote{
Note that in this standard theory, one inverts $L_\sigma -\mu I$ rather
than $\cD^2-1$.},
we find that there exists a closed ball $B(\lambda_*, \eps)$ of radius $\eps$, centered at $\lambda_*$, containing
no other eigenvalues of $L_\sigma$.

Consider now an increasing sequence of eigenspaces $\{H_J\}_{j\in\mathbb{N}}$ of $L^2_{\rm per}([0,X])$
such that $\lim_JH_J=L^2_{\rm per}([0,X])$ and $H_*\subset H_J$ for all $J\in\mathbb{N}$.
For each $J$, let $\{r_k\}_{k=1}^J$ be an orthonormal basis of $H_J$ and let $R_J=(r_1,\ldots,r_J)$.
Then we can define the
finite-dimensional approximants
\be\label{factn}
\delta_J(\sigma,\lambda):= {\det}_{2}\left( R^*_J (\partial_x^2-1)^{-1}(L_\sigma-\lambda I)R_J\right).
\ee
Since $D_\sigma $ does not vanish on $\partial B(\lambda_*,\eps)$,
by the correspondence in location of roots and eigenvalues established above, and since
$\delta_J$ converges locally uniformly in $\lambda$ to $D_\sigma$ by \eqref{compare},
Rouch\'e's Theorem implies that there exists a $J^*\in\NM$ sufficiently large such that for $J>J^*$ 
the winding number of $D_\sigma $ around $\partial B(\lambda_*,\eps)$
is equal to the winding number of $\delta_J$ around the same ball.

Finally, fixing $J_0>J^*$ and noticing that $L_\sigma R_{J_0}=R_{J_0} M_{\sigma,J_0}$, where $M_{\sigma,J_0}$
is an $J_0\times J_0$ matrix representation of $L_\sigma$ on the finite-dimensional
invariant subspace $H_{J_0}$, we find from \eqref{factn} that there exists a constant $C\neq 0$ such that
\[
\delta_{J_0}(\sigma,\lambda)={\det}_{2}\left(R_{J_0}^*(\partial_x^2-1)^{-1}R_{J_0}  (M_{\sigma,J_0}-\lambda I)\right)=
C{\det}_{2}\left(M_{J_0}-\lambda I\right),
\]
and hence we see that $\delta_{J_0}$
is a nonvanishing multiple of the characteristic polynomial
of $M_{\sigma,J_0}$.  Here, we are using the fact that $R_{J_0}^*(\partial_x^2-1)^{-1}R_{J_0}$
is positive definite, by positive symmetric definiteness of
$(\partial_x^2-1)^{-1}$.
It follows that $\delta_{J_0}$ 
has a zero at $\lambda_*$ of precisely the algebraic
multiplicity of $\lambda_*$ as an eigenvalue of $L_\sigma$.
Thus, we conclude that the multiplicity of $\lambda_*$ as a root of $D_\sigma$ is
equal to the winding number of $\delta_{J_0}(\cdot,\sigma)$ about the ball $\partial B(\lambda_*,\eps)$,
which in turn is equal to the algebraic multiplicity of $\lambda_*$
as an eigenvalue of $L_\sigma$, completing the proof.
\end{proof}

\br\label{ges}
The truncated winding-number argument for agreement of multiplicity
to our knowledge is new, and seems of general use in similar situations.
It would be interesting to prove this also in a different way by establishing
a direct correspondence between the Fredholm determinant and the standard
periodic Evans function construction of Gardner \cite{G}, as done in the
solitary-wave case in \cite{GLM1, GLMZ2, GM} and in the periodic
Schr\"odinger case in \cite[Sect.\  4]{GM}.
This would give at the same time an alternative proof of Gardner's fundamental
result of agreement in location and multiplicity of roots of
the standard periodic Evans function with eigenvalues of $L_\sigma$,
through the result of Theorem \ref{evansthm}.
\er

\subsection{Convergence of Hill's method}\label{s:hill}

Next, we
use the machinery developed in the previous section to give a
proof of the convergence of Hill's method.  
In order to precisely describe Hill's method,
notice that by taking the Fourier transform,
we may express \eqref{eig} equivalently as
the infinite-dimensional matrix system
\[
(\D^2 +\D\A_1 + \A_0- \lambda I)\U=0,
\]
where for each $m=0,1$ and $j,k\in\ZM$,
\be\label{mats}
\D_{jk}=\delta_j^k  ij,
\quad
[\A_m]_{jk}= \widehat{A_m}(j-k),
\quad\textrm{ and }
\U_j= \widehat U(j),
\quad
\ee
where $\hat f(k)$ denotes the discrete Fourier tranform of $f$ evaluated at Fourier frequency $k$ and,
as elsewhere, $i=\sqrt{-1}$.
Hill's method then consists of fixing $J\in\NM$ and truncating the above
infinite-dimensional matrix system
at wave number $J$, that is,
considering the $(2J+1)$-dimensional minor $|(j,k)|\le J$,
and computing the eigenvalues of the finite-dimensional matrix 
\begin{equation}\label{approxL}
L_{\sigma,J}:=\D_J^2+\D_J\A_{1,J}+\A_{0,J},
\end{equation}
where $\D_J$ and $\A_{m,J}$ denote the $(2J+1)$-dimensional matrices resulting from truncating the matrices $\D$ and $\A_m$
to frequencies $|(j,k)|\leq J$,
to obtain approximate eigenvalues for $L_\sigma$.
Notice this can be done 
quite efficiently by applying modern numerical linear
algebra techniques. 


\br\label{divrmk}
In
applications, one may of course
 encounter operators $L$ that are not in divergence
form \eqref{e:L}.
In this case, we point out that there is no effect in changing from
nondivergence to divergence form except that we increase the regularity
requirement on $A_1$ from $L^2$ to $H^1$.
Indeed,, we may change
from one form to the other using the Leibnitz rule
$\A_1\D-\D\A_1=  (\A_1)'$, where
$$
(\A_1)'_{jk}= i(j-k)\A_1(j-k)=
(\widehat {A_{1,x}})(j-k),
$$
and noting that, 
since $\D$ is diagonal, this operation is respected by truncation.
Thus, there is indeed no loss of generality in our representation of
operators in divergence form, as it does not affect the result of
Hill's method.
\er

Following the construction of the generalized periodic Evans function \eqref{evans2},
we may rewrite the truncated eigenvalue equation
\be\label{trunceig}
 \left(L_{\sigma,J}-\lambda I\right)\cU=0
\ee
as
\be\label{truncfred}
 (I+ \cK_J)\cU=0,
\ee
where $\cK_J=\cK_{1,J}+\cK_{2,J}$ is the truncation of the
Fourier representation $\cK=\cK_1+\cK_2$ of operator $K$ to frequencies $|(j,k)|\leq J$, that is,
\be\label{trun}
 \cK_{1,J}=\D_J(\D_J^2-I)^{-1}\A_{1,J}
 \quad\textrm{ and }\quad
 \cK_{2,J}=(\D_J^2-I)^{-1}(\A_{0,J}+1-\lambda).
\ee
Continuing to follow the above construction of $D_\sigma$, we now define the truncated
periodic Evans function as
\be\label{truncevans}
D_{\sigma,J}(\lambda):={\det}_{2}(I-\cK_J)
\ee
and notice that we have the following preliminary result.

\bl\label{trunccorr}
The zeros of $D_{\sigma,J}$ correspond
in location and multiplicity with those of $L_{\sigma,J}$.
\el

\begin{proof}
This is immediate by the
nonsingularity of $(\cD_J^2-I)^{-1}$ and properties
of the (usual, finite-dimensional) characteristic polynomial,
together with the observation that
\[
{\det}_2(I-\cK_J)=
{\det}_2(\cD_J^2-I)^{-1} {\det}_2(\D_J^2 +\D_J\A_{1,J} + \A_{0,J}-\lambda I).
\]
\end{proof}

With this construction in hand, we now state the main result of this section.

\bt\label{convthm}
For $A_j\in L^2_{\rm per}([0,X])$, the sequence of determinants $D_{\sigma,J}$
converges to $D_\sigma$ as $J\to \infty$ uniformly in $\sigma$ and locally uniformly in $\lambda$.
\et

\begin{proof}
This convergence result follows from the proof of Theorem \ref{evansthm}.
Indeed, 
noting
that $D_{\sigma,J}$ is exactly such a
sequence of approximate determinants, corresponding here 
to the ascending sequence of sinusoidal functions of integer wave number,
by which the generalized periodic Evans function $D_\sigma$ was
defined in \eqref{evans2}, 
we find by our definition of the 2-modified Fredholm determinant that
$D_{\sigma,J}\to D_\sigma$ pointwise in $\lambda$
as $J\to \infty$ for each fixed $\sigma\in[0,2\pi)$. 
Moreover, recalling that the rate of convergence is determined by
the difference between truncated operator $\cK_J$ and $\cK$ in
$\cB_2$
norm, and noting that we have uniformly bounded
$\cB_2$
estimates on each entry of $\cK_J$, 
we find that this convergence is uniform in $\sigma$ and locally uniform in $\lambda$. 
\end{proof}

From Theorem \ref{convthm} we immediately have convergence of Hill's method, as described in the introduction.
For completeness, we state this result in the following corollary.

%
%
\bc\label{Lconvthm}
For $A_j\in L^2_{\rm per}([0,X])$,
 the eigenvalues of $L_{\sigma,J}$ defined in \eqref{approxL}
approach the eigenvalues of $L_\sigma$ in location and multiplicity
as $J\to \infty$, uniformly on $|\lambda|\le R$, $\sigma\in [0,2\pi]$,
for any $R$ such that $\partial B(0,R)$ contains no eigenvalues of
$L_\sigma$.
\ec

\begin{proof}
This is immediate from Theorem \ref{evansthm}, Lemma \ref{trunccorr}, and Theorem \ref{convthm},
along with basic 
properties of uniformly convergent analytic functions.
\end{proof}

\subsection{Rates of Convergence}

Next, we address the issue of the rates of convergence of $D_{\sigma,J}$ to $D_\sigma$ and of the
approximate spectra to the exact spectra.  Assuming slightly more regularity on the function $A_1$
in \eqref{e:L}, we have the following easy convergence result.

\bt\label{Lconvrate}
For $A_j\in H^1_{\rm per}([0,X])$ and each fixed $R>0$, there exists a constant $C=C(R)>0$ such that
for each fixed $|\lambda|\leq R$
\[
|D_{\sigma,J}(\lambda)- D_\sigma(\lambda)|\le CJ^{-1/2}.
\]
In particular, this estimate is locally uniform in $\lambda$ and uniform in 
$\sigma$.
\et

\begin{proof}
The rate of convergence is bounded by
$\|\cK_J-\cK\|_{\cB_2}$
from which we readily obtain the result using the Cauchy-Schwarz
estimate
$$
\sum_{|j|\ge J}|\widehat{A^m}(j)|^2\le \sum_{|j|\ge J}|j|^{-2}
\sum_{|j|\ge J}|j|^{2}|\widehat{A^m}(j)|^2\le (C/J)\|A^m\|_{H^1([0,X])}
$$
for each $m\in\{0,1\}$.
For details, see the very similar estimates in the proof of
Theorem 4.9, \cite{GLZ}.
\end{proof}

Notice that Theorem \ref{Lconvrate} does not imply
a rate of convergence of the roots of $D_{\sigma,J}$ to the
roots of $D_\sigma$, or, equivalently, the eigenvalues of $L_{\sigma,J}$ to
the eigenvalues of $L_\sigma$.
%
Indeed, the above convergence result
is, with or without rate information, essentially an abstract one.
For, though we find convergence the of analytic functions
$D_{\sigma,J}$ to $D_\sigma$,
we don't obtain rates of convergence of their zeros without more structural
information about $D_\sigma$ itself.  In particular, we can not conclude
convergence rates of the approximate spectra to the true eigenvalues of $L_\sigma$
using only the knowledge of the 
eigenvalues of $L_{\sigma,J}$ computed in the course of Hill's method.
This suggests the idea of computing the approximate Evans function
$D_{\sigma,J}$ directly, instead of using it as a purely analytical tool,
an idea that would be interesting for future investigation.  Though
in principle slower due to the need for multiple evaluations of eigenvalues,
this computation is better conditioned, so there might perhaps be some
counterbalancing advantages to this approach, besides the possibility
already mentioned to obtain a posteriori estimates on the error bounds
for eigenvalue approximations.
We leave this as an interesting topic for further investigation,
related to the larger question of relative advantages
of standard periodic Evans function (as in \cite{G})
vs. Hill's computations.

\section[Generalizations]{Generalizations}
\label{s:comp}

Here, we briefly discuss various generalizations of the theory developed in
Section \ref{s:simple}.

\subsection{Operators with nontrivial principal coefficient}\label{s:princ}

Consider now a system of the more general form
\be\label{e:Ln}
L_\sigma= \partial_x^2 A_2 + \partial_x A_1(\sigma,x) + A_0(\sigma,x),
\ee
where $A_2$ is symmetric positive definite, satisfying $ A_2(x)\ge C$
for some $C>0$, uniformly on $x\in [0,X]$.
Define as usual $\cA_2$ to be the infinite-dimensional
matrix representation of $A_2$ under Fourier transform; that is,
$\cA_{2,jk}=\widehat{A_2}(j-k)$.
Then clearly $\cA_2$ is symmetric and, by Parseval's identity,
satisfies $ \cA_2 \ge C$ when considered as a quadratic
form on $\ell^2(\NM)$.
As a consequence, the $J^{\rm th}$ truncation $\cA_{2,J}$, as a principal
minor of a positive definite symmetric matrix, must also be positive definite
and satisfy the same bound $ \cA_{2,J}\ge C$.

In particular, $\cA_2$ is invertible with
$$
 \cA_2^{-1} \ge 1/C,
\quad
 \cA_{2,J}^{-1} \ge 1/C.
$$

\bl\label{multlem}
$\|AB\|_{\cB_2}\le |A|_{L^2}\|B\|_{\cB_2}$,
where $|\cdot|_{L^2}$ denotes $L^2([0,X])$ operator norm.
\el

\begin{proof}
Straightforward from the definition of $\|\cdot\|_{\cB_2}$.
\end{proof}

\bc\label{c:compose}
For $A_j\in L^2_{\rm per}([0,X])$ and $A_2$ symmetric positive definite with $A_2(x)\ge C$,
the operator $\cM:=\cA_2^{-1}\cK$ is Hilbert-Schmidt where $\cK=\cK_1+\cK_2$ is defined
as in \eqref{trun}.
\ec

In this case, following the notation of Corollary \ref{c:compose}, we
define the generalized Evans function as $D_\sigma(\lambda):=\det_2 (I-\cM)$,
noting that the eigenvalue problem may be written equivalently
as $(I-\cM)\cU=0$.
The associated series of Fredholm approximants is
$ D_{\sigma,J}(\lambda):=\det_2(I-\cM_J)$, with $D_{\sigma,J}(\lambda)
\to D_\sigma(\lambda)$ uniformly as $J\to \infty$, just as before,
and zeros of $D_{\sigma}$ corresponding
in location and multiplicity with eigenvalues of $L_\sigma$.
However, the corresponding object obtained by Hill's method is not
the truncated Fredholm determinant $D_{\sigma,J}$ defined above, but rather the modified version
\be\label{hillapprox}
\check{D}_{\sigma,J}(\lambda):= {\det}_2( I- \cA_{2,J}^{-1} \cK_J),
\ee
and it is this function whose zeros correspond with the eigenvalues
of the Hill approximant operator $L_{\sigma,J}$.

To verify convergence of Hill's method in this case then, it is sufficient to show that
\be \label{key}
\| \cM_J- \cA_{2,J}^{-1} \cK_J\|_{\cB_2}
=\| (\cA_2^{-1}\cK)_J- \cA_{2,J}^{-1} \cK_J\|_{\cB_2}
\to 0
\ee
as $J\to \infty$.  Indeed, with this convergence result in hand
we may conclude 
by \eqref{compare} that $\lim_{J\to\infty}|\check{D}_{\sigma,J}-D_{\sigma,J}|= 0$, 
and
thus $\check D_{\sigma,J} \to D_\sigma$ as $J\to \infty$, yielding
the convergence result as before.

\bt
For operators of the form \eqref{e:Ln}, Hill's method converges in location and multiplicity
provided that $A_j\in L^2_{\rm per}([0,X])$.
\et

\begin{proof}
We sketch the proof of \eqref{key}.
By boundedness of $\|A_2\|_{L^2([0,X])}$, we may truncate $\widehat{A_2}$
at wave number $M$ to obtain an $M$-banded infinite-dimensional diagonal matrix
centered around zero-frequency approximating $\cA_2$ to arbitrarily small order in the $\ell^2(\NM)$
operator norm.  Hence, for purposes of this argument, we may assume
without loss of generality that $\cA_2$ is $M$-banded diagonal operator centered about zero-frequency.
Furthermore, noting that since $\widehat{\cA_2^{-1}}$ is bounded in $L^2(\RM)$, for $J\in\NM$ sufficiently large 
the columns of $\cA_2^{-1}$ corresponding to frequencies $|j|\leq J-M$ are small off the principal $2J+1-M$ minor
and hence
a brief calculation reveals
that
$$
(\cA_2^{-1})_J \cA_{2,J}=
\bp
E_{M}& 0 & 0\\
0 & I_{2J-2M} & 0\\
0 & 0 &F_{M}\\
\ep,
$$
where $E_M$ and $F_M$ are $M\times M$ matrices that are invertible
by invertibility of
$(\cA_2^{-1})_J \cA_J$, a property of principal minors of positive-definite
symmetric matrices.
By a further left-multiplication by the block-diagonal matrix
$$
\bp
E_{M}^{-1}& 0 & 0\\
0 & I_{2J-2M} & 0\\
0 & 0 &F_{M}^{-1}\\
\ep
$$
we obtain
$I_{2J+1}$,
demonstrating that
$(\cA_{2,J})^{-1}$ agrees with $(\cA^{-1})_{2,J}$ on the central
$2J-2M+1$ dimensional minor.
Recalling that $\|\cK-\cK_J\|_{\cB_2}\to 0$ as $J\to 0$ by \eqref{e:kconverge},
we thus obtain by a straightforward calculation
$$
\| (\cA_2^{-1}\cK)_J- \cA_{2,J}^{-1} \cK_J\|_{\cB_2}
\sim \| (\cA_2^{-1}\cK_J)_J- \cA_{2,J}^{-1} \cK_J\|_{\cB_2}
\to 0,
$$
completing the proof by \eqref{compare}
\end{proof}

\subsection{Composite and Higher-order operators}

The reader may easily verify that all of the arguments of Sections
\ref{s:simple} and \ref{s:princ} carry over to the case when the operator \eqref{e:L1}
is replaced by a general periodic-coefficient operator
$$
L=\partial_x^m a_m(x) + \partial_x^{m-1}a_{m-1}(x)+\dots +a_0(x)
$$
where $a_j\in L^2_{\rm per}([0,X])$ and where the principal coefficient $a_m$ symmetric positive definite.
Indeed, the analysis parallels that of previous sections except that one must
substitute for $(\partial_x^2-1)$ 
everywhere the positive definite symmetric
Fourier multiplier
\[
|\partial_x^2-1|^{m/2}= \cF^{-1}(|j|^2+1)^{m/2}\cF,
\]
where $j$ denotes the Fourier wave number
and $\cF$ denotes Fourier transform.  With these substitutions, our previous arguments
immediately yield convergence of Hill's method in this case as well. 


Furthermore, it is straightforward to verify that all of the analysis
of Sections \ref{s:simple} and \ref{s:princ} extends readily to the case of operators
of ``composite'' type
\[
L=\bp
\partial_x^{m_1}a^1_{m_1} +\dots\\
\vdots\\
\partial_x^{m_n}a^n_{m_n} +\dots\\
\ep,
\]
with $a^j_{k}\in L^2_{\rm per}([0,X])$ and $a^j_{m_j}$ symmetric positive definite
for each suitable choice of indices:
that is,
still assuming $L$ is a nondegenerate ordinary differential operator in some sense.

\br\label{compapp}
It is the above observation that applies to the numerics in \cite{BJNRZ1,BJNRZ2}, where
the authors use Hill's method to numerically analyze the spectrum of the linearized
St. Venant equations 
\begin{align*}
\lambda \tau-c\tau' - u'&= 0,\\
\lambda u-cu'
-(\bar \tau^{_-3}(F^{-1}- 2\nu \bar u_x)\tau)' &=
-(s+1)\bar \tau^s\bar u^r \tau
- r\bar \tau^{s+1}\bar u^{r-1} u
+\nu (\bar \tau^{-2}u')' 
\end{align*}
about a given periodic or homoclinic orbit $(\bar u,\bar \tau)$,
where $r$, $s$, $F$, and $\nu$ are physical parameters in the problem and $\lambda$
is the corresponding spectral parameter.
\er

\subsection{Operators with general coefficients}
Our results are completely general in the scalar case,
applying to all nondegenerate operators.
However, they are restricted in the system case by the
condition that the principal coefficient(s) be symmetric positive
definite.  Whether this condition may be relaxed is an interesting operator-theoretic
question regarding properties of Toeplitz matrices.

Specifically, the property that we need to carry out Hill's method
(and indeed, to complete our entire convergence analysis) is
that the minor $\cA_{2,J}$ of a Toeplitz matrix $[\cA_{2}]_{mn}=\widehat{A_2}(k-n)$
be invertible for $J$ sufficiently large.
The question is what properties of $A_2(x)$ are sufficient
to guarantee this: in particular, is uniform invertibility enough?
Alternatively, what are sufficient conditions on $\widehat{A_2}?$
This seems an interesting problem for further investigation.

%

\medskip
{\bf Acknowledgement.} Thanks to Bernard Deconink for pointing
out the references \cite{CuD,CDKK,DK}.

\end{document}